\documentclass{amsart}
\usepackage[utf8]{inputenc}
\usepackage{amssymb}
\usepackage{amsmath}
\usepackage{array}
\usepackage{hyperref}
\usepackage{xypic}
\usepackage{multirow}
\usepackage{xcolor}

\usepackage{tikz}
\usepackage{tikz-cd}
\usetikzlibrary{shapes.geometric, arrows}

\newcommand{\Lk}{\mbox{$\mathfrak k$}}

\newcommand{\Lh}{\mbox{$\mathfrak h$}}

\newcommand\fieldsetc{\mathbb}

\newcommand{\R}{\fieldsetc{R}}

\newtheorem{thm}{Theorem}[section]
\newtheorem{cor}[thm]{Corollary}
\newtheorem{prop}[thm]{Proposition}
\newtheorem{lem}[thm]{Lemma}

\theoremstyle{remark}
\newtheorem{rmk}[thm]{Remark}

\newcommand{\Pf}{{\em Proof}. }
\newcommand{\EPf}{\hfill$\square$}
\title{Inhomogeneous almost symmetric submanifolds}

\author{Claudio Gorodski and Carlos Olmos}\thanks{The first named author
acknowledges partial support from FAPESP Thematic Grant
22/16097-2 and CNPq Fellowship 304252/2021-2, Brazil. The second named
author acknowledges partial support by FaMAF,
Universidad Nacional de Córdoba and CIEM-CONICET, Argentina.}

\dedicatory{To our dear colleague and friend Paolo Piccione on his 60th
  birthday.}

\subjclass{53C40, 53C42, 53C35}

\keywords{Almost symmetric submanifolds; cohomogeneity one
submanifolds; extrinsic de Rham's theorem}

\begin{document}

\maketitle

\begin{abstract}
  We completely describe inhomogeneous properly embedded almost symmetric
  submanifolds of Euclidean space as certain unions
of parallel symmetric submanifolds of the ambient Euclidean space.

\end{abstract}

\section{Introduction}

Recall that a submanifold of a Euclidean space is called \emph{(extrinsically) symmetric}
if it is reflectionally symmetric with respect to any affine normal
space. Symmetric submanifolds were introduced and classified
by Ferus~\cite{Fe}, who locally characterized them by the
condition that the second fundamental form be parallel (see also~\cite{St,ET}).
In particular Ferus showed that symmetric submanifolds are homogeneous and
arise as special orbits of the isotropy representation of a symmetric space,
or \emph{s-representation}, for short. Symmetric submanifolds
split extrinsically, and the irreducible ones are classified
by the work of Kobayashi and Nagano~\cite{KN}. These
submanifolds enjoy a rich interplay between extrinsic geometry and representation theory, as the geometry is closely tied to the structure of the isotropy representation.

Almost symmetric submanifolds of Euclidean space form a closely related
class. They are interesting because, despite relaxing full symmetry, they
retain strong rigidity properties. 
An \emph{(extrinsic) almost symmetry} of a Euclidean submanifold
$M$ at a point $p$
in $M$ is an involutive ambient isometry that preserves $M$,
fixes the affine normal space at~$p$ pointwise,
and whose fixed point set in the tangent
space has dimension one. A Euclidean submanifold is called \emph{(extrinsically)
almost symmetric} if it admits an almost symmetry at every point.
We introduced this class of submanifolds in~\cite{GO}, proved that
they are of cohomogeneity at most one, and classified the homogeneous ones.
Herein we complete that discussion by giving a complete description of the inhomogeneous case.

Our main result reads as follows.

\begin{thm}\label{main}
  Consider an irreducible s-representation of a compact connected Lie group 
  $K_i$ on $V_i$ with a fixed
  non-trivial extrinsic symmetric  orbit $K_i\cdot v_i$, for $i=1, \ldots , r$.
  Put $V = V_0\oplus V_1\oplus \cdots \oplus V_r$, where $V_0$ is a 
possibly trivial Euclidean space,  and consider the diagonal action of 
$K=K_1\times \cdots \times K_r$ on~$V$. We consider a 
connected properly embedded $1$-dimensional
  submanifold $L$ of $V_0 \oplus
    \mathbb{R} v_1\oplus \cdots \oplus \mathbb{R}v_r$ through
$v=v_1+\cdots+v_r$ such that  
the orthogonal projection of $L$ to $\mathbb{R}v_i$ is nonconstant
and does not meet
the origin, for all $i=1, \ldots , r$.
Moreover, we suppose that 
the orthogonal projection of $L$ to~$V_0$  is full in $V_0$ . Then:
\begin{enumerate}
\item[1.] $KL$ is a full irreducible properly embedded
  almost symmetric submanifold of $V$.
\item[2.] Fix a nonzero vector $w\in V_0$. Denote by $w^\perp$ the 
hyperplane of $V_0 \oplus
    \mathbb{R} v_1\oplus \cdots \oplus \mathbb{R}v_r$ which is orthogonal to $w$. 
Suppose that 
$L$ is not contained in $w^\perp$
and that it is 
invariant under the reflection in that hyperplane. 
Then $(SO(k)\times K)L$ is a full irreducible properly embedded
almost symmetric submanifold of 
$\R^{k-1}\oplus V$, where $SO(k)$ acts by the standard representation on
$\R^{k-1}\oplus\R w$ and trivially on its orthogonal complement, and $k\geq2$. 
    \item[3.] Fix nonzero orthogonal vectors $w$, $w'\in V_0$.
Denote by $w^\perp$ (resp.~$w^{\prime\perp}$) the 
hyperplane of $V_0 \oplus
    \mathbb{R} v_1\oplus \cdots \oplus \mathbb{R}v_r$ which is orthogonal to $w$
(resp.~$w'$).  Suppose that 
$L$ is neither contained in $w^\perp$ nor in~$w^{\prime\perp}$,   
and its invariant under the reflection in either hyperplane. 
Then $(SO(k)\times SO(k')\times  K)L$ is a full irreducible properly
embedded almost symmetric submanifold of~$\R^{k-1}\oplus\R^{k'-1}\oplus V$, where $SO(k)$ (resp.~$SO(k')$) acts by the standard representation on
$\R^{k-1}\oplus\R w$ (resp.~$\R^{k'-1}\oplus\R w')$ and trivially on its orthogonal 
complement, and $k$, $k'\geq2$.  
\end{enumerate}    

Further, the submanifolds constructed above are inhomogeneous, unless they are a round sphere. 
Conversely, every full irreducible 
inhomogeneous properly embedded almost symmetric submanifold of Euclidean space is constructed as above. 
\end{thm}

\section{Preliminary material}

\begin{lem}\label{lem:23b}
Consider an oriented connected complete one-dimensional
submanifold $L$ of $\mathbb{R}^n$ ($n \geq 2$) which is 
invariant under the reflection $r$ on a codimension one subspace
$W$ of $\R^n$, but it is not contained in $W$. 
Denote by $e$ be the positive unit tangent vector field along~$L$.  Then:
\begin{enumerate}
    \item $L\cap W$ is non-empty. 
    \item If $rq=q$ for some $q\in L$, then $re_q=-e_q$ 
and hence $r|_L$ is the geodesic symmetry of $L$ at $q$.
    \item $L\cap W$ has at most two elements;
equivalently, $r$ has at most two fixed points. 
    \item If $L\cap W$ consists of a single point, 
then $L$ is non-compact. 
    \item If $L\cap W=\{p,q\}$, then: $L$ is compact; $p$, $q$ 
are $L$-antipodal points; and $e_p= -e_q$.
    \item Suppose further that $L$ is invariant under  
the reflection $\bar r$ on another
 codimension one subspace
$\bar W$ of $\R^n$ such that $r$ and $\bar r$ commute,
but $L$ is not contained in $\bar W$. 
Then $L$ is compact. Moreover, if $L\cap W=\{p,q\}$,
$L\cap \bar W=\{\bar p,\bar q\}$, then the four points 
$p$, $\bar p$, $q$, $\bar q$ 
are mutually distinct and equally spaced along~$L$.
\end{enumerate}
\end{lem}

\Pf Denote by $W^\pm$ the open half-spaces defined by $W$. 
Since $rL=L\not\subset W$, we have that
$L\cap W^+$ and $L\cap W^-$ are non-empty disjoint subsets 
of $L$. Now $(a)$ follows from connectedness of~$L$.

From $rL=L$ and $rq=q$, we see that $dr_qe_q=re_q=\pm e_q$.
But $r$ is not the identity along $L$, so it must be that $re_q=-e_q$. 
This proves~(b). 

Suppose $p\in L$ is fixed by $r$.  Then $r|_L$ is the
 geodesic symmetry at $p$. Now $r$ has another fix point along $L$
if and only if $L$ is 
compact, in which case it must be the antipodal point to $p$, that is, the 
farthest point to $p$. In this case $L\cap W=\{p,q\}$ and we can take 
an arc-length parameterization of $L$ given
by~$c:[0,\tau]\to \mathbb{R}^n$ where 
$c(0)=p$, $c(\tau/2)=q$ and $c'(0)=e_p$. Also, by relabeling $W^\pm$ 
if necessary, we 
may assume that  $c(t)\in\mathbb{W}^+$ for $t\in(0,\tau/2)$.
Owing to~(b), $r|_L$ is the geodesic symmetry at $p$, so
$r(c(\tfrac{\tau}{2}-t)\bigr) = c(\tfrac{\tau}{2}+t)$. This yields
that  $c(t)\in\mathbb{W}^-$ for $t\in(\tau/2,\tau)$. 
But $r|_L$ is also the geodesic symmetry at~$q$, so 
$c'(\tau/2)\perp W$. It follows that $e_q=c'(\tau/2)=-c'(0)=-e_p$.
This proves (c), (d) and~(e). 

If $p$ and $\bar p$ are arbitrary points of~$L$
fixed by $r$ and $\bar r$, respectively, then 
$p\neq \bar p$ because $e_p\perp W$, $e_{\bar p}\perp \bar W$ and 
$W\neq \bar W$. In particular, $\bar r p\neq p$. 
But $r$ and $\bar r$ commute, so $q:=\bar r p\in L$ is a 
fix point of $r$. By~(e), $L$ is compact. Since the geodesic
symmetry at any of $p$, $\bar p$, $q$, $\bar q$ must permute these
four points, they must be equally spaced along~$L$.
This completes the proof of~(f) and the proof of the lemma. 
\EPf

\begin{cor} 
Let $W_1,\ldots, W_k$ be mutually orthogonal linear 
hyperplanes in $\R^n$ for $3\leq k\leq n$. 
Consider a connected complete one-dimensional
submanifold $L$ of $\mathbb{R}^n$ 
which is invariant under the reflection $r_i$ on $W_i$
for $i=1,\ldots,k$. Then 
$L$ is contained in~$k-3$ of the hyperplanes $W_1,\ldots,W_k$.
\end{cor}

\begin{lem}\label{preimage}
Let $C$ be a connected properly embedded $1$-submanifold of $\R^n$ 
which is invariant under the reflection $r$ on a hyperplane $W$ of $\R^n$
and which is not contained in~$W$. 
Then there exists a $r$-invariant 
neighborhood $\Omega$ of $C$ and a $r$-invariant submersion
$f:\Omega\to\R^{n-1}$ such that $C=f^{-1}(0)$.

If, in addition, 
$C$ is invariant under  the reflection $\bar r$ on a second hyperplane
of $\R^n$ which is orthogonal to the first one,
and $C$ is not contained in $\bar W$, then $\Omega$ and $f$ can be taken 
$\bar r$-invariant.  
\end{lem}

\Pf  Owing to Lemma~\ref{lem:23b},
there is $p\in C$ such that $rp=p$. Let $\xi_1,\ldots,\xi_{n-1}$
be an orthonormal basis of~$\nu_pC=W$
and extend it to a parallel orthonormal frame 
$\tilde\xi_1,\ldots,\tilde\xi_{n-1}$ of normal vector fields
along $C$. Note that this frame 
is well-defined, even if $C$ meets the hyperplane at a second point, 
since indeed each $\tilde\xi_i$ is $r$-related to itself. This
shows that the normal bundle of $C$ in $\R^n$ is globally flat,
and in particular trivial.  

Set
\[ \Gamma(t,x_1,\ldots,x_{n-1})=c(t)+\sum_{i=1}^{n-1}x_i\tilde\xi_i, \]
where $c:J\to\R^n$ is a parameterization of $C$ and $J=S^1$ or $\R$. 
It is clear that $\Gamma$ has full rank along $x_1=\ldots=x_{n-1}=0$. 
Since $C$ is properly embedded, there exists a smooth function 
$\epsilon:J\to\R$ such that $\Gamma$ is a diffeomorphism from 
\[ B_\epsilon(0) =\{(t,x_1,\dots,x_{n-1})| \sum_i x_i^2<\epsilon(t),\ t\in J\} \]
onto its image $B_\epsilon(C)$. Note that the function $\epsilon$
can also be taken $r$-invariant, so that $B_\epsilon(C)$ 
becomes a $r$-invariant open neighborhood of $C$. 
Now
\[ f=\pi\circ\Gamma^{-1}:B_\epsilon(C)\to\R^{n-1}, \]
where $\pi:\R^n\to\R^{n-1}$ is projection onto the last $n-1$ coordinates, 
is a submersion such that $f^{-1}(0)=C$. Since each $\tilde\xi_i$ is $r$-related 
to itself, it holds that $f$ is $r$-invariant. For the last assertion, 
we can take $\epsilon$ to be $\bar r$-invariant. Since each $\tilde\xi_i$ is
also $\bar r$-related to itself, the result follows.  \EPf

\begin{lem}\label{preimage2}
  Let $C$ be as in the first part of
  Lemma~\ref{preimage} and assume that $W$ is the hyperplane 
$e_1^\perp$, where $e_1$ is the first vector in the canonical basis of $\R^n$.
Then $(SO(k)\times \{I\})C$ is a properly embedded submanifold
of $\R^n\oplus\R^{k-1}$, where $SO(k)$ acts as usual 
on $\R e_1\oplus\R^{k-1}$, and trivially on its orthogonal complement
in $\R^n\oplus\R^{k-1}$. 
If, in addition, $C$ is invariant under the reflection $\bar r$
on $\bar W=e_2^\perp$,
then $(SO(k)\times SO(k')\times \{I\})C$ is a properly embedded submanifold
of $\R^n\oplus\R^{k-1}\oplus\R^{k'-1}$,
where $SO(k)$ (resp.~$SO(k')$) acts as usual 
on $\R e_1\oplus\R^{k-1}$ ($\R e_2\oplus\R^{k'-1}$), and trivially on its orthogonal complement in $\R^n\oplus\R^{k-1}\oplus\R^{k'-1}$.
\end{lem}

\Pf The $r$-invariance of~$f$ means that 
\[ f(-x_1,x_2,\ldots,x_n)=f(x_1,x_2.\ldots,x_n) \]
for $(x_1,\ldots,x_n)\in\Omega$.  
Then there is a smooth function $h$ such that 
\[ f(x_1,x_2,\ldots,x_n)= h(x_1^2,x_2,\ldots,x_n) \]
(see e.~g.~\cite{Sch}). 
We just define $\tilde f:\tilde\Omega\to\R^{n-1}$ by
\[ \tilde f(x_1,\ldots,x_n,y_1,\ldots,y_{k-1})=h(x_1^2+y_1^2+\cdots+y_{k-1}^2,x_2,\ldots,x_n), \]
where $\tilde\Omega$ is a suitable open subset of $\R^{n+k-1}$,
and note that $\tilde f$ is smooth and
$(SO(k)\times \{I\})C=\tilde f^{-1}(0)$ is a regular level set.

If $C$ is also $\bar r$-invariant, we can take $h$ satisfying
\[ f(x_1,x_2,x_3,\ldots,x_n)= h(x_1^2,x_2^2,x_3,\ldots,x_n), \]
and we can define
$\tilde f:\tilde\Omega\to\R^{n-1}$ by
\begin{eqnarray*}
  \lefteqn{\tilde f(x_1,\ldots,x_n,y_1,\ldots,y_{k-1},z_1,\ldots,z_{k'-1})}\\
&&=h(x_1^2+y_1^2+\cdots+y_{k-1}^2,
x_2^2+z_1^2+\cdots+z_{k'-1}^2,
x_3,\ldots,x_n),
\end{eqnarray*}
where $\tilde\Omega$ is a suitable open subset of
\[\R^n\times \R^{k-1}\times\R^{k'-1}= \R^{n+k+k'-2},\]
and the result about $(SO(k)\times SO(k')\times \{I\})C$ follows as above. \EPf

\section{Characterization of inhomogeneous almost symmetric submanifolds}\label{charact}

Let $M$ be a properly embedded inhomogeneous
almost symmetric submanifold of Euclidean space~$V$.  Denote
by $K$ the closure of the
group of isometries of $M$ generated by the almost symmetries of $M$,
By one of the main results in~\cite{GO},
we know that $K$ with cohomogeneity $1$ on $M$.
For a $K$-regular point $p\in M$, the almost symmetry $\sigma_p$ is unique
and can be regarded as a symmetry of the Euclidean submanifold
$S=K^0p$. It is also true that $K^0$ is the (extrinsic) transvection group of
$S$, namely, the closure of the group generated by products of two
symmetries of~$S$. 

It follows from the work of Ferus (see~\cite[Theorem~2.8.8]{BCO})
that $S$ splits extrinsically as the product of an Euclidean space
and a compact symmetric submanifold, and the group $K^0$ splits
off a group of translations along a subspace of $V$. Its orbits
define a parallel distribution of $V$, whose restriction along $M$
is tangent to~$M$. Then it is a standard application of Moore´s lemma
to deduce that $M$ splits an Euclidean factor. 
At this juncture
we assume that $M$ is full and irreducible. This implies that there
is no Euclidean factor and $S$ is just a compact symmetric submanifold. 
Now there are splittings $K^0=K_1\times\cdots\times K_r$,
$V=V_0\oplus V_1\oplus\cdots\oplus V_r$ such that $V_0$ is a
trivial representations and $K_i$ acts as an irreducible
s-representation on $V_i$ and trivially on $V_j$ for $j\neq i$,
for all $i=1,\ldots,r$.

We decompose $p=p_0+p_1+\cdots+p_r$ where $p_i\in V_i$. It is
known that each irreducible symmetric submanifold $K_ip_i$ for $i\neq0$
is a most singular orbit, so its orbit type is isolated in the sphere of
$V_i$ of radius $||p_i||$. 

Let $\xi$ be a unit normal vector to~$S$ in $M$ at~$p$,
and consider the normal
geodesic $\gamma_\xi$ to $S$ starting at $p$. The fixed point set
of $(d\sigma_p)_p$ on $T_pM$ is $\R\xi$. It follows that
the image
$C$ of $\gamma_\xi$ is a closed totally geodesic submanifold of $M$, 
and the almost symmetry
$\sigma_{\gamma_\xi(t)}$ coincides with $\sigma_p$ whenever $\gamma_\xi(t)$ is a
$K^0$-regular point.

\begin{lem}\label{lem:C}
With the notation above, 
\[ C\subset V_0\oplus\R p_1\oplus\cdots\oplus\R p_r. \]
Moreover,
\begin{equation}\label{comps}
 \gamma_\xi(t)=c_0(t)+\lambda_1(t)p_1+\cdots+\lambda_r(t)p_r, \
\end{equation}
where $c_0$ is a smooth curve in $V_0$ through $p_0$, and $\lambda_i$ is a 
smooth function with $\lambda_i(0)=1$ for~$i=1,\ldots,r$.
Further, $\lambda_i$ never vanishes unless $(K_i,V_i)=(SO(k),\R^k)$;
in the latter case $\lambda_i(t_0)=0$ for some $t_0\neq0$ precisely
when $\gamma_\xi(t_0)$ is a $K^0$-singular orbit.
\end{lem}

\Pf We know that $\sigma_p$ is an almost symmetry of $M$ at
$\gamma_\xi(t)\in C$ for all $t\in\R$. By $\sigma_p$-equivariance
of the second fundamental form $\alpha$ of $M$ in $V$, we get that 
$\alpha(\gamma_\xi'(t),\mathcal H_t)=0$ for all~$t$, where 
$\mathcal H_t=\gamma_\xi'(t)^\perp\cap T_{\gamma_\xi(t)}M$. 
Further, $\mathcal H_t$ is parallel along $\gamma_\xi$, since 
the latter is a geodesic of $M$. It follows that $\mathcal H_t$ is a 
constant subspace $\mathcal H_0=T_pS$ of $V$. We deduce that 
$\gamma_\xi'(t)\in\nu_pS$ for all $t$. Of course $S$ is a principal
orbit, so $C$ is poinwise fixed by the isotropy group $(K^0)_p$, which 
yields that $(K^0)_p$ fixes $\gamma_\xi'(t)$ for all $t$. But $S$ is an
extrinsic symmetric submanifold of $V$, from where we know that the 
fixed point set of $(K^0)_p$ on $\nu_pS$ is exactly the fixed point set 
$\nu^0_pS$ of the restricted normal holonomy of $S$ 
at $p$ in $V$. So $C\subset p+\nu^0_pS$.  
Since irreducible symmetric submanifolds are most singular orbits, 
the $(K^0)_p$-fixed point set in $\nu_pS$ is
\[ V_0\oplus\R p_1\oplus\cdots\oplus\R p_r. \]
Note that $p$ belongs to this subspace, which proves the first and 
second assertions. 

Finally, it is clear that $\lambda_i(t_0)=0$ for some $i$
if and only if $K^0\gamma_\xi(t_0)$ is a
singular orbit. Suppose this is the case and write $q=\gamma_\xi(t_0)$,
$N=K^0q$. By~\cite{GO}, the
elements of $K^0$ that act trivially on $N$ form
a connected normal subgroup $H$ that is isomorphic to $SO(k)$, where
$k\geq2$ is the codimension of $N$ in $M$. Since $N$ is totally geodesic
in $M$, its normal bundle in $M$ is parallel. Also, the $H$-equivariance
of the second fundamental
form $\alpha$ of $M$ in $V$ implies that $\alpha(T_xN,\nu_xN\cap T_xM)=0$
for all $x\in N$. It follows that the normal bundle of $N$ in $M$
is parallel in $V$, so it is a constant subspace $W$ of $V$.
Since $H$ is generated
by the composition of any two almost symmetries of $M$ at $q$, $H$ also
acts trivially on $\nu_xM$. But $H$ is a normal
subgroup of $K^0$, so $K^0$ preserves the fixed point set $W^\perp\subset V$
of $H$. 
Now it is clear that $K_i=H\cong SO(k)$ and $V_i=W\cong\R^k$. \EPf

\begin{rmk}
We keep the notation in the proof of Lemma~\ref{lem:C}, 
and recall that the subgroup $H$ of $K^0$ is independent of the point $x\in N$.
Since a cohomogeneity one action has at most two singular orbits, say
$N$ and $N'$, this
argument shows that there are at most two factors $SO(k)$, $SO(k')$
of $K^0$ associated to $N$, $N'$. Note that even if there are two
singular orbits, their normal factors can coincide (e.g.~surfaces of
revolution). 
\end{rmk}

\begin{lem}
  With the notation above, suppose $\lambda_i(t_0)=0$ so that
  $q=\gamma_\xi(t_0)$ is a
$K^0$-singular point and $K_i$ the associated $SO(k)$-normal factor
of $K^0$ as in Lemma~\ref{lem:C}. Then $C$ is invariant under the orthogonal 
reflection on the hyperplane $p_i^\perp$
of $V_0\oplus\R p_1\oplus\cdots\oplus\R p_r$ orthogonal to $p_i$. 
\end{lem}

\Pf Note first that $\gamma_\xi'(t_0)=\lambda_i'(t_0)p_i\in V_i$,
where of course $\lambda_i'(t_0)\neq0$. 
This is because $\gamma_\xi$ is orthogonal to every orbit it meets, 
and the normal space to $K^0q$ in $M$ is precisely $V_i$.   
This implies that $p_i\in V_i$. 
Let $g\in SO(V_i)\subset (K^0)_q$ map $\gamma_\xi'(t_0)$ to its opposite.
Then $gp_i=-p_i$.  
It also is clear that $gC=C$, and we claim that $g$ is the reflection
on $p_i^\perp\subset V_0\oplus\R p_1\oplus\cdots\oplus\R p_r=\nu_p^0S$.
Indeed, 
\[ \Lk_ip_i=\Lk_ip\subset\Lk p=T_pS, \]
and $g$ is the identity on
\[ V_i^\perp=(\Lk_ip_i\oplus\R p_i)^\perp\supset\nu_pS\cap p_i^\perp\supset
\nu^0_pS\cap p_i^\perp=p_i^\perp, \]
as wished. \EPf

\medskip

Recall that $C=\gamma_\xi(\R)$ meets all the $K^0$-orbits. It follows that
$M=K^0C$. Further, $C\subset\nu_p^0S$, the fixed point set
of the restricted normal holonomy. Therefore
\[ M=\bigcup_t S_t, \]
where $S_t=K^0\gamma_\xi(t)$ is a parallel submanifold to $S$ in $V$.

We refer to equation~(\ref{comps}) and remark that 
the fullness of the component curve
$c_0$ in $V_0$ is equivalent to the fullness of $M$ in $V$. In fact,
each s-orbit $K_ip_i$ is full in $V_i$. For each $i\neq0$
it is also clear that the
orthogonal projection of $C$ to $\R p_i$ is nonconstant, equivalently
each $\lambda_i$ is nonconstant, for otherwise $M$ would be an extrinsic
product
\[ (K^0/K_i)\cdot(p-p_i)\times K_i\cdot p_i, \]
and we have assumed $M$ irreducible. This completes the proof of the last
assertion in the statement of the theorem. 

\section{Construction of inhomogeneous almost symmetric submanifolds}

In this section we prove 1, 2 and 3 in the statement of the theorem.

\begin{lem}
Let $M$ be equal to
$KL$ (resp.~$(SO(k)\times K)L$, $(SO(k)\times SO(k')\times K)L$)
be as in the statement of Theorem~\ref{main}. Then  $M$ is a 
properly embedded (resp.~properly embedded, compact)
almost symmetric submanifold of $V$ 
(resp.~$\R^{k-1}\oplus V$, $\R^{k-1}\oplus \R^{k'-1}\oplus V$).
\end{lem}

\Pf Consider $M=KL$ and fix $v\in L$. 
Then the map $K/K_v\times L\to V$, $(kK_v,w)\mapsto kw$ 
is well-defined (since the isotorpy group of $K$ is constant along $L$) 
and a proper embedding (since $L$ is properly embedded and $K$ is compact). 
This shows that $KL$ is a properly embedded submanifold of $V$. 

Let $\sigma_v$ be the symmetry of $S=Kv$ as a submanifold of $V$. Since 
\[ L\subset V_0\oplus\R v_1\oplus\cdots\oplus\R v_r=\nu_v^0S=v+\nu_v^0S
\subset v+\nu_vS, \]
$L$ is pointwise fixed by $\sigma_v$. Then 
\[ \sigma_vKL=(\sigma_vK\sigma_v^{-1})\sigma_vL=KL, \]
which implies that $\sigma_v$ is an almost symmetry of $M$ at~$v$. 
Note that $\sigma_v$ is an almost symmetry of $KL$ at any point of $L$. 

Next we consider $M=(SO(k)\times K)L$. We note first that 
due to Lemma~\ref{preimage2}
$N=(SO(k)\times \{I\})L$ is a properly embedded submanifold 
of $\R^{k-1}\oplus V_0\oplus\R v_1\oplus\cdots\oplus\R v_r$. 
Then we can emply the reasoning above to see that 
$M=KN$ is a properly embedded submanifold of $V$.
Note that $(SO(k)\times K)v$ is a symmetric submanifold of 
$V$, namely, the product of a $(k-1)$-sphere with $Kv$, which has
codimension one in $M$, if $v$ is a $SO(k)\times K$-regular point of $M$.
 As in the previous case, the symmetry of $(SO(k)\times K)v$ at $v$
is an almost symmetry of $M$. 

The last case we consider is  $M=(SO(k)\times SO(k')\times K)L$.
Note that $N=(SO(k)\times SO(k')\times \{I\})L$
is a compact submanifold 
of $\R^{k-1}\oplus\R^{k'-1}\oplus V_0\oplus\R v_1\oplus\cdots\oplus\R v_r$
again by Lemmata~\ref{lem:23b} and~\ref{preimage2}.
Similarly to the above, one shows
that $M=KN$ is a
compact submanifold of~$V$, and $(SO(k)\times SO(k')\times K)v$ is a
symmetric submanifold of $V$, which has codimension one in $M$, if
$v$ is a $(SO(k)\times SO(k')\times K$-regular point of~$M$;
the symmetry of this submanifold is an almost symmetry of~$M$.
\hfill\mbox{ }\hfill \EPf

\begin{lem}
  Let $M$ be equal to
$KL$ (resp.~$(SO(k)\times K)L$, $(SO(k)\times SO(k')\times K)L$)
be as in the statement of Theorem~\ref{main}. Then  $M$ is an
(extrinsically) irreducible submanifold of~$V$ 
(resp.~$\R^{k-1}\oplus V$, $\R^{k-1}\oplus \R^{k'-1}\oplus V$).
\end{lem}

\Pf Suppose $M=KL=M'\times M''$ is an extrinsic product of submanifolds
of Euclidean space. For a $K^0$-regular point $q=(q',q'')\in M$
with $q'\in M'$, $q\in M''$, consider the almost symmetry $\sigma_q$.
By Theorem~\ref{thm:uniq-extr-factors}, $\sigma_q$ must permute the
extrinsic factors, but $\sigma_q$ fixes the affine
normal space of~$M$ at~$q$ pointwise, so indeed $\sigma_qM'=M'$ and $\sigma_qM''=M''$.
It follows that the $+1$-eigenspace of $d(\sigma_q)_q$ must be contained
in either $T_{q'}M'$ or $T_{q''}M''$. Say we are in he first case, so that
$\sigma_q$ restricts to an (extrinsic)
almost symmetry of $M'\times\{q''\}$, and it
restricts to an (extrinsic) symmetry of $\{q'\}\times M''$.
At points nearby  $q$, the set of fixed points of the unique
almost symmetry changes
smoothly, as it coincides with the normal spaces to the $K^0$-orbits.
Again by Theorem~\ref{thm:uniq-extr-factors}, $\sigma_{(q',x)}$
restricts to a symmetry of $\{q'\}\times M''$ for $x$ near $q''$.
We claim that this implies
that $\{q'\}\times M''$ is a symmetric (homogeneous) 
submanifold of the ambient Eucldiean
space. In fact, the local one-parameter groups of global transvections
of  $\{q'\}\times M''$ along geodesics emanating from $q$ can be extended
to a transitive group of isometries, since $\{q'\}\times M''$ is complete.
Now, due to the identity component of the isometry group of $M$ being the product of the identity components of the isometry groups of $M'$ and $M''$,
the homogeneity of  $\{q'\}\times M''$  implies that it is contained
in $S=K^0q$. Denote by $\mathcal D$ the parallel distribution of $M$
associated to the factor $M''$. Since $K^0$ is connected,
it must preserve $\mathcal D$. In particular, we see from the
above inclusion that $\mathcal D$ is tangent to $S$ along $S$. 
Now it is a standard application of Moore's lemma to deduce that
there is an extrinsic product $S=S'\times M''$.
If we take $q=v$ and write $v=(v',v'')$ where $v'\in M'$, $v''\in M''$,
then it follows from that decomposition that $\{v'\}\times M''=K''v$,
where $K''$ is the product of some factors of $K^0$. Since $L$ is tangent
to $M'\times\{v''\}$, it projects constantly to $\{v'\}\times M''$ and
hence to some $\R v_i$ with $i\neq0$, which is a contradiction to
our assumption.

In the cases of $M=(SO(k)\times K)L$
(resp.~$M=(SO(k)\times SO(k')\times K)L$), we argue in the same way
just by replacing the group $K$ by
$SO(k)\times K$ (resp.~$SO(k)\times SO(k')\times K$), and in the end
we note that $L$ cannot project constantly to $\R^k$ (resp.~
$\R^k\oplus\R^{k'}$) since it passes through $0$ in $\R^k$
(resp.~in~$\R^k$ and in~$\R^{k'}$). \EPf

\medskip

The reasoning in the end of section~\ref{charact} yields that
$M$ is a full submanifold of Euclidean space, and this completes
the proof of parts 1, 2 and 3 of the theorem. 

\section{Inhomogeneity}

Finally, we show that the submanifolds constructed in parts 1, 2 and 3
of the theorem are inhomogeneous, unless they are a round sphere.
For that purpose, we shall use Simon's theory of holonomy systems
to study enlargements
of s-representations. 

Let a compact connected Lie group 
$K$ act as an s-representation on $V$, up to a (possibly non-existent)
trivial component $V_0$. There is an
algebraic curvature tensor $R$ associated
to the representation of $K$ on $V$, namely, $R$ is $K$-invariant 
and the Lie algebra $\Lk$ of $K$ is spanned by $R_{x,y}$ for $x$, $y\in V$.  
Further, $R$ can be assumed to be positive semidefinite, since Cartan duality
preserves the equivalence class of the isotropy representation of a
symmetric space. 

\begin{lem}\label{h}
With the above definitions,
let $H$ be a closed connected subgroup of $SO(V)$ that contains $K$. 
Then $H=H'\times H''$ and $V=V^{H''}\oplus W$, 
where $V^{H''}\subset V_0$ is the fixed point set of $H''$,
$H''$ acts as an s-representation on $W$, and $H'$ acts trivially 
on $W$. In particular, if $V^{H''}$ is nonzero or if $H''$ does not 
act irreducibly, then every $H$-orbit in $V$ is either reducible or 
not full. 
\end{lem}

\Pf Define $H(R)$ to be the linear span of the $h\cdot R$ for $h\in H$,
and $\Lh^R$ to be the subspace of $\Lh$ spanned by $S_{x,y}$ for 
$S\in H(R)$ and $x$, $y\in V$. Then $\Lh^R$ is an ideal of $\Lh$.
Denote by $H^R$ the corresponding connected (normal) subgroup of $H$. 
Write 
\begin{equation}\label{w}
V=W_0\oplus W_1\oplus \cdots\oplus W_s,
\end{equation}
$H^R$-irreducible decomposition, where $W_0$ is the fixed point set
of $H^R$.  The symmetries of algebraic curvature tensors yield that 
every $S\in H(R)$ satisfies
$S_{W_0,W_i}=0$ for all $i$, $S_{W_i,W_j}=0$ for $i\neq j$, and 
$S_{W_i,W_i}W_i\subset W_i$ for all $i$. This implies that  
the linear span $\Lh_i^R$ of $S_{x,y}$ for $x$, $y\in W_i$, is 
an ideal of $\Lh^R$ that acts irreducibly on $W_i$, trivially on $W_j$ for 
$j\neq i$, $\Lh_0^R=0$, and  $\Lh^R=\Lh^R_1\oplus\cdots\oplus\Lh^R_s$
(direct sum of ideals). It follows that the decomposition~(\ref{w})
is unique, up to the order of the factors. Since $H$ normalizes
$H^R$, it preserves each $W_i$,  and it 
is irreducible on $W_i$ for $i\neq0$. 

Next, for each $i\neq0$, we choose a nonzero $R^i\in\Lh^R_i$.
Since $R^i$ is semidefinite positive as a curvature endomorphism,
its scalar curvature $sc(R^i)$ is positive, and therefore 
$[W_i,R^i,H^R_i]$ is an irreducible holonomy system with $sc(R^i)\neq0$, 
where $H^R_i$ is the connected (normal) subgroup of $H^R$ with 
Lie algebra $\Lh^R_i$. By~\cite[Theorem~5]{simons}, $H^R_i=H|_{W_i}$ 
acts as an s-representation on $W_i$. Defining $W=W_1\oplus\cdots\oplus W_s$,
this yields that $H|_W=H^R|_W$. Now
given $h\in H$, there is $h^\circ\in H^R$ such that $h(h^\circ)^{-1}$ 
is the identity
on $W$. Writing $h=(h(h^\circ)^{-1})h^\circ$, we see that $H=H'\times H''$,
where $H''=H^R$, as wished. \EPf

\begin{prop}
Let a compact connected Lie group 
$K$ act as an s-representation on $V$, up to a (possibly non-existent)
trivial component. Suppose an orbit $M=Kv$ has codimension~$1$
in $Hv$, where $H$ is a closed connected subgroup of $SO(V)$ containing $K$,
and where $Hv$ is full and irreducible in $V$. Then $Hv$ is the round sphere
of radius $||v||$. Moreover, if $M$ is a symmetric submanifold of $V$,
then $M$ must be the product of round spheres or a non-full round
sphere of codimension one in $Hv$.
\end{prop}

\Pf Due to Lemma~\ref{h},
$H=H^R$ acts as an irreducible s-representation on~$V$. Fix an
algebraic curvature tensor $R\neq0$ associated to $K$. Then
$[V,R,H]$ is an irreducible holonomy system. It cannot be symmetric, for that
would spell that $\Lh=\Lh^R$ is spanned by $R_{x,y}$ for $x$, $y\in V$,
and therefore equals $\Lk$. Owing to~\cite[Theorem~4]{simons}, we
know that $H$ is transitive on the unit sphere of $V$.

For the last statement, recall that a homogeneous hypersurface of a sphere
$S^{n+1}$ is isoparametric, and irreducible isoparametric hypersurfaces
of the sphere are not symmetric submanifolds, so $M$ must be a reducible
isoparametric hypersurface of $S^n$. By the classification
$M$ is either $S^{m_1}\times S^{m_2}$ with $m_1+m_2=n$, or
an umbilic $S^n$ in $S^{n+1}$. \hfil\mbox{ }\hfil\EPf 

\section{Appendix: Uniqueness of extrinsic splittings of isometric
immersions}

The well-known Moore's lemma states that an isometric immersion $f$ of a 
Riemannian product manifold decomposes into an extrinsic product 
of isometric immersions of the factors if and only if the second
fundamental form of $f$ decomposes accordingly. There is also a 
uniqueness result, for which we provide a proof now, 
since we could not locate one in the literature.

\begin{thm}\label{thm:uniq-extr-factors}
Let $f:M\to\R^n$ be a full isometric immersion of a symply-connected
complete Riemannian manifold. Suppose that $f$ decomposes as 
an extrinsic product of isometric immersions
\[ M=M_0\times M_1\times \cdots \times M_r,\quad f=f_0\times\cdots\times f_r, \]
where $M_0=\R^{k_0}$ is an Euclidean space ($k_0\geq0$)
and $f_i:M_i\to\R^{k_i}$  are irreducible isometric 
immersions with $k_i\geq2$ 
for $i=1,\ldots,r$. Then this decomposition is unique, up to a permutation of 
the non-Euclidean factors. 
\end{thm}

\Pf Let $\mathcal D_0, \dots, \mathcal D_r$ be the parallel distributions 
in $M$ associated with its product decomposition. 
Observe that the distribution $\mathcal D^{dR}_0$, associated with the flat de 
Rham factor $\mathbb R^k$ of $M$, intersects non-trivially $\mathcal D_i$ 
for some $i \geq 1$ if $k_0 < k$. Fix a point $p$ in $M$. 
Let $v\in T_p\mathcal D^{dR}_0(p)$ and let $\tilde v$ its associated parallel vector field on $M$. Observe that any parallel vector field on $M$ must be tangent to its flat de Rham factor. It is standard to show that that $v\in \mathcal D_0$ if and only if $df(\tilde v_x)$ is a constant vector of $\mathbb R^n$, i.e., it
is independent of $x$. This shows that $\mathcal D_0$ does not depend on the decomposition of the immersion. The factor  $M_0= \mathbb R^{k_0}$ is called the extrinsic flat factor of the isometric immersion $f$. 
Thus we may assume that $k_0=0$.

Suppose there is another product decomposition 
$M = \tilde M_1 \times \cdots \times \tilde M_s$
with associated isometric immersions 
$\tilde f_j : \tilde M_j \to \R^{\tilde k_j}$. 
Then the tangent and normal bundles of $M$ decompose in two different ways, 
corresponding to the two products of immersions. Namely, we can write
$$ TM_1 \oplus \dots \oplus TM_r = TM = T\tilde M_1 \oplus \dots \oplus T\tilde M_s \, ,$$
$$ \nu M_1 \oplus \dots \oplus \nu M_r = \nu M =  \tilde \nu M_1 \oplus \dots \oplus \tilde \nu M_s \, ,$$
with the following properties:
\begin{itemize}
    \item[(i)] The tangent bundles $TM_i$ and $T\tilde M_j$ are invariant under all shape operators~$A_\xi$. Equivalently, by the Gauss formula, 
    \[
    \alpha(TM_i, TM_i^\perp) = 0 = \alpha(T\tilde M_j, T\tilde M_j^\perp).
    \]
    
    \item[(ii)] 
    \[
    A_{\xi_i} TM_{i'} = 0 = A_{\tilde \xi_j} T\tilde M_{j'} \quad \text{for } i\neq i',\ j\neq j',
    \]
    where $\xi_i$ and $\tilde \xi_j$ are sections of $\nu M_i$ and $\nu \tilde M_j$, respectively.
\end{itemize}

Define
\[
\mathcal F = \left\{ \tau_c A_{\xi_x} \tau_c^{-1} : \xi_x \in \nu_x M,\ c \text{ an arbitrary curve from } x \text{ to } p, \ x \in M \right\},
\]
where $\tau _c$ denotes parallel transport in $TM$ along $c$. 
For $i = 1, \dots, r$, $j = 1, \dots, s$, 
similarly define $\mathcal F_i$ and $\tilde{\mathcal F_j}$ by respectively
restricting $\xi_x$ to be a vector in $\nu_xM_i$, $\nu_x\tilde M_j$.
Then 
\begin{equation*}\label{eq:ew1} \mathcal F_1 \oplus \dots \oplus \mathcal F_r = \mathcal F = 
\tilde{\mathcal F}_1 \oplus \dots \oplus \tilde{\mathcal F}_s,
\end{equation*}
\begin{equation}\label{eq:ew2}    \mathcal F_i (TM_{i'}) = 0 = \tilde{\mathcal{F}_j} (T\tilde M_{j'}) \quad \text{for } i\neq i',\ j\neq j', 
\end{equation}
and
\begin{equation*}\label{eq:ew3}
\begin{split}
\mathcal F(TM_i) &= \mathcal F_i(TM_i) \subset TM_i, \\
\mathcal F(T\tilde M_j) &= \tilde{\mathcal F}_j(T\tilde M_j) \subset T\tilde M_j.
\end{split}
\end{equation*}

Observe that $\ker\mathcal F$, the common kernel of all the elements of $\mathcal F$, is  invariant under the holonomy of group of $M$ at $p$. Therefore
 it defines a parallel distribution $\mathcal D$ of $M$ such that $\alpha(\mathcal D, TM) = 0$, or, equivalently, 
$\alpha(\mathcal D, \mathcal D^\perp) = 0$ and 
$\alpha(\mathcal D, \mathcal D) = 0$. By Moore's lemma  
$M=M'  \times M'' $, where $M'$ and $M''$ are the factors of $M$ 
associated to $\mathcal D$ and $\mathcal D^\perp$, respectively,
and $f=f'\times f''$ where $f'$ and $f''$ 
 are isometric immersions into Euclidean spaces. 
Moreover, since $\alpha(\mathcal D, \mathcal D) = 0$,  $f'$
is totally geodesic. Then $M'$ is an extrinsic flat factor of $f$,
contrary to our assumption, unless $\mathcal D =0$. 
This shows that  $\ker \mathcal F =0$, and implies that 
$\mathcal F_i(T_pM_i)=\mathcal F(T_pM_i)\neq 0$ and 
 $\tilde{\mathcal F}_i(T_p\tilde M_i)=\mathcal F(T_p\tilde M_i)\neq 0$, for all $i=1, \dots , r$, $j = 1 , \dots , s$.

We deduce that for each $i=1, \dots , r$ there is 
$j=1, \cdots , s$ such that 
\[
0\neq \mathcal F_i (T_p\tilde M_j)=\mathcal F(T_p\tilde M_j)
\subset T_p\tilde M_j.
\]
From (\ref{eq:ew2}), 
the non-trivial subspace $\mathcal F_i (T_p\tilde M_j)$ lies in 
$T_pM_i\cap T_p\tilde M_j$. This intersection is invariant under both $\mathcal F$ and the holonomy of $M$. Hence it defines a parallel distribution, 
which moreover is invariant under all shape operators of $M$ 
due to the definition of $\mathcal F$. 
Owing to Moore’s lemma, it corresponds to a non-trivial extrinsic factor. 
Such a factor must be contained in both $M_i$ and $\tilde M_j$, which are irreducible extrinsic factors. Therefore 
$M_i=\tilde M_j$. Arguing by induction, this yields that $r=s$ and 
that $M_i=\tilde M_{\sigma(i)}$, for some permutation $\sigma$ 
of $1,\ldots, r$. \EPf

\providecommand{\bysame}{\leavevmode\hbox to3em{\hrulefill}\thinspace}
\providecommand{\MR}{\relax\ifhmode\unskip\space\fi MR }
\providecommand{\MRhref}[2]{%
  \href{http://www.ams.org/mathscinet-getitem?mr=#1}{#2}
}
\providecommand{\href}[2]{#2}

%\bibliographystyle{amsalpha}
%\bibliography{ref}

\end{document}